\documentclass[11pt,a4paper]{amsart}

\usepackage[utf8]{inputenc}
\usepackage{amsmath,amssymb,amsthm}
\usepackage{color,graphicx,ulem,url,tikz}
\RequirePackage{hyperref}

\newtheorem{theorem}{Theorem}
\newtheorem{definition}[theorem]{Definition}
\newtheorem{example}{Example}

\newcommand*\circled[1]{\tikz[baseline=(char.base)]
{\node[shape=circle,draw,inner sep=1pt](char){#1};}}

\begin{document}


\title[Hospital management in the COVID-19 emergency]{Hospital management in the COVID-19 emergency: Abelian Sandpile paradigm\\ and beyond}

\author[R. Martucci]{Roberta Martucci}
\address[Roberta Martucci]{Dipartimento di Matematica G. Castelnuovo, Sapienza Universit\`a di Roma, piazzale Aldo Moro 2 - 00185 Roma (Italy)}
\email{martucci.1701560@studenti.uniroma1.it}

\author[C. Mascia]{Corrado Mascia}
\address[Corrado Mascia]{Dipartimento di Matematica G. Castelnuovo, Sapienza Universit\`a di Roma, piazzale Aldo Moro 2 - 00185 Roma (Italy)}
\email{corrado.mascia@uniroma1.it}

\author[C. Simeoni]{Chiara Simeoni}
\address[Chiara Simeoni]{Laboratoire de Math\'ematiques J.A. Dieudonn\'e CNRS UMR 7351, Universit\'e C\^ote D'Azur, Parc Valrose - 06108 Nice Cedex 2 (France)}
\email{chiara.simeoni@univ-cotedazur.fr}

\author[F. Tassi]{Filippo Tassi}
\address[Filippo Tassi]{Dipartimento di Matematica G. Castelnuovo, Sapienza Universit\`a di Roma, piazzale Aldo Moro 2 - 00185 Roma (Italy)}
\email{tassi.1912485@studenti.uniroma1.it}


\begin{abstract}
In this article, we propose a mathematical model --based on a cellular automaton-- for the redistribution of patients within a network of hospitals with limited available resources, in order to reduce the risks of a local/global collapse of the healthcare system. We attempt at developing a conceptual tool to support making rational decisions relevant to the optimisation of the allocation of patients into accessible medical facilities. The strategy is based on a version of the Abelian Sandpile model for the Self-Organised Criticality, with the idea of testing the paradigm for the management of patients among the COVID-19 hospitals in Italian regions. In particular, we compare the novel proposal to the standard management of connections between hospitals, showing a number of advantages at a local and global level, by means of a reliable indicator function introduced for measuring the effectiveness of the allocation strategies.
\end{abstract}



\maketitle


\section{Introduction}
\label{section:intro}
\smallskip

\begin{figure}[htb]
\includegraphics[width=0.85\textwidth]{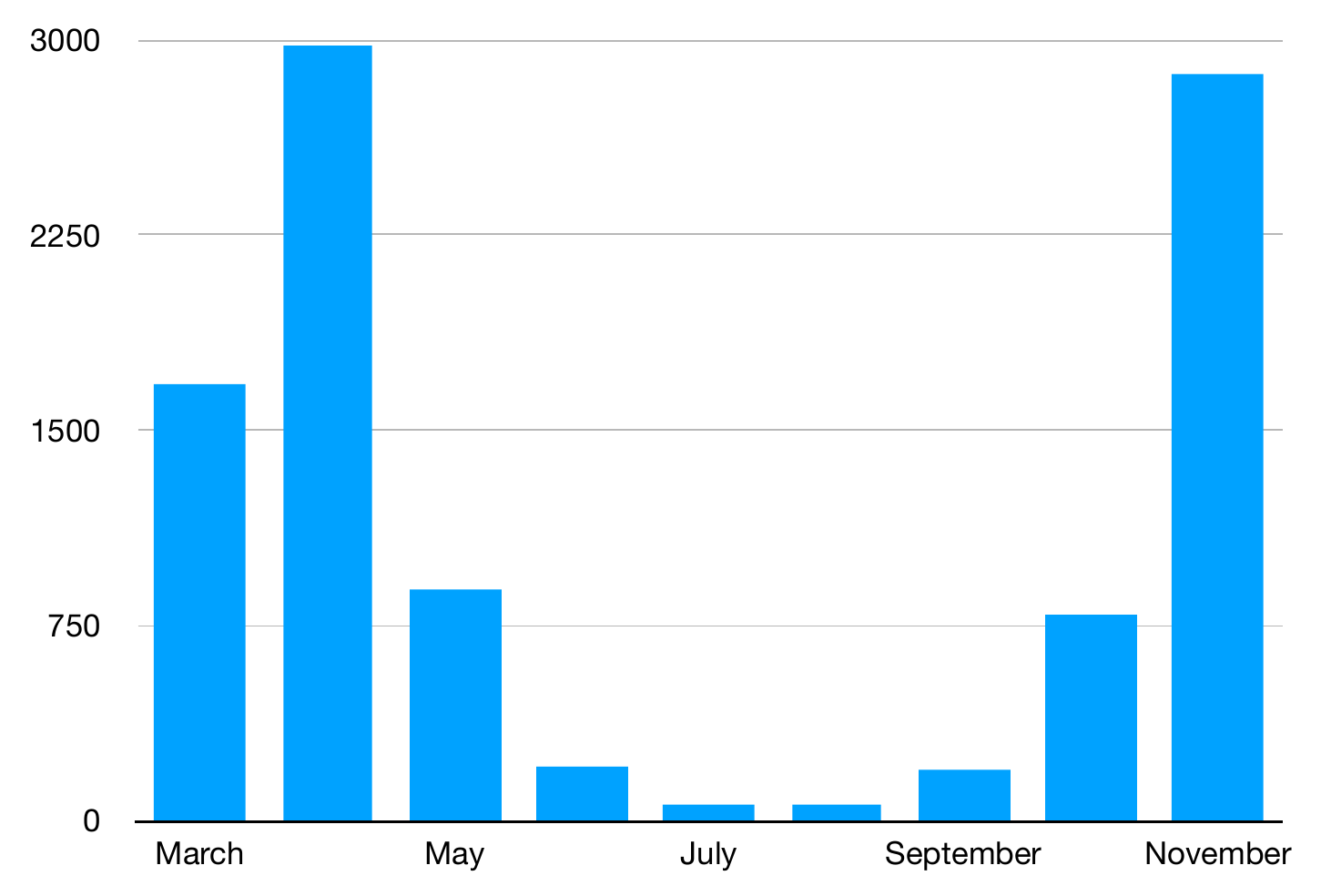}
\caption{\small Occupancy level (aggregate data) of intensive care units in the Italian healthcare system.}
\label{figure:TerapieIntensive}
\end{figure}

Since the beginning of 2020, people from all parts of the world are struggling with the COVID-19 coronavirus pandemic, which has led to a hundred million of infections and over 2 million deaths worldwide~\cite{ECDC2021,WHO2021}. The health emergency is not over yet and it will take months before being under control, possibly becoming endemic among most of the world population~\cite{LBA2021}.\\
The current data suggest the paramount importance of acting promptly, and especially for trying to reorganize the medical facilities in order to avoid saturation of the critical care capacity in the various territories (refer to Figure~\ref{figure:TerapieIntensive} for the evolution of employment of the intensive care units in Italy). Indeed, several countries have been experiencing recurring epidemic waves, with a consequent worsening of the social and economic conditions, thus increasing dramatically the pressure on the healthcare systems.

\begin{figure}[htb]
\includegraphics[width=0.95\textwidth]{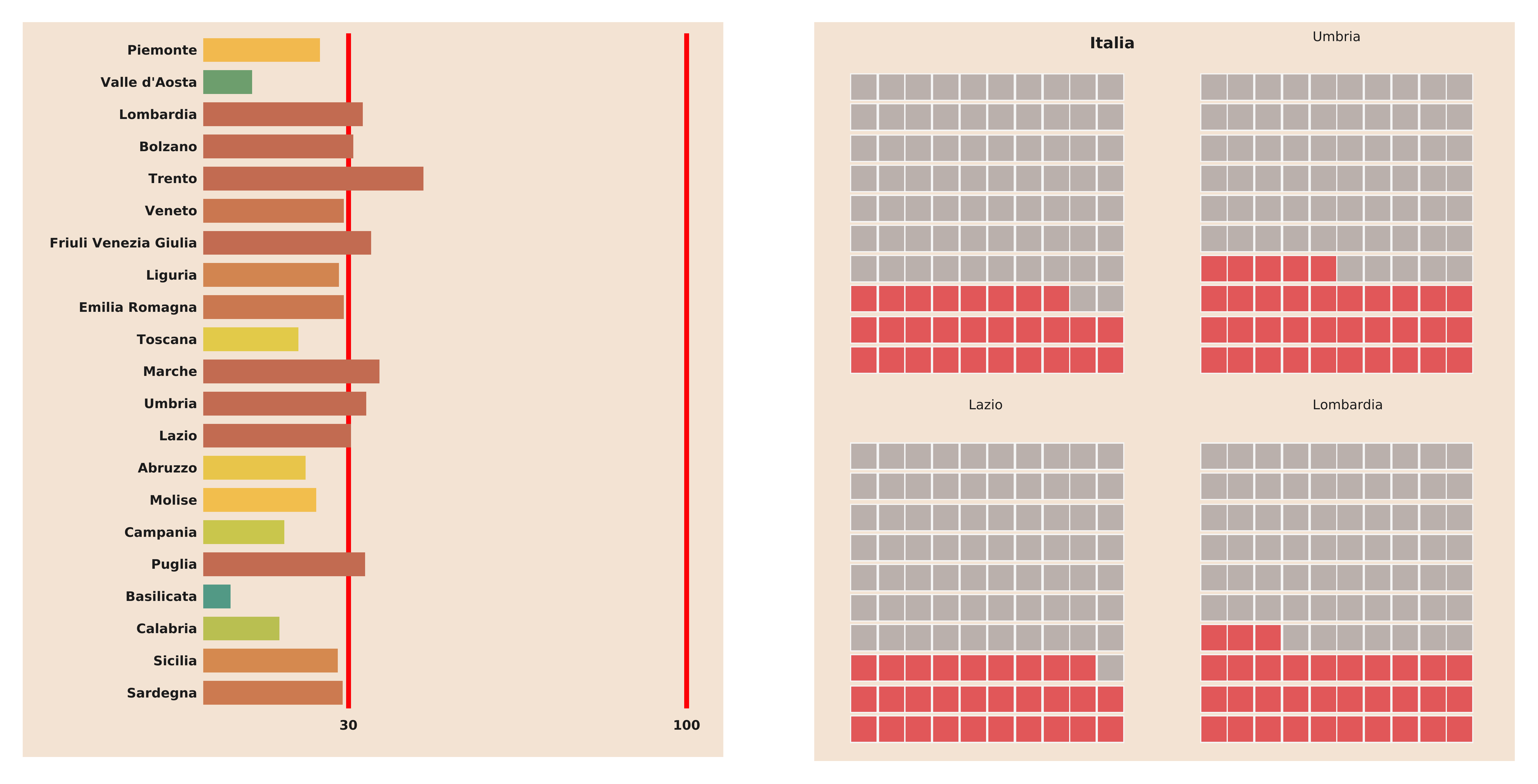}
\caption{\small Saturation of intensive care units in the Italian regions (updated to January 23, 2021).
Source: Il Sole 24 Ore -- elaboration on data from the Italian Civil Protection Department}
\label{figure:saturazione}
\end{figure}

The emergence of infectious diseases has to be considered as an inherent property of human and animal populations, which cannot be generally avoided and/or foreseen: before the present COVID-19 emergency, SARS (2002-2004) and influenza A/H1N1 (2009-2010) are other recent notable incidents~\cite{WIK2002,WIK2009}.\\
Most epidemics are characterised by outbreaks in localized geographical areas, and there may be no pharmaceutical solutions already available to safeguard the people's health when the disease begins to spread largely around. Moreover, the development of effective vaccines typically takes a rather long time, and therefore temporary alternative strategies must be implemented (quarantine, travel restrictions, activities closing, social distancing, face masks obligation, ...) However, the demand for specific healthcare facilities such as hospitals consistently increases, reaching and sometimes exceeding critical levels. For example, in the case of COVID-19, the saturation of intensive care units depends on the amount of patients requiring mechanical support for ventilation and proper equipment, so that any healthcare system can be overwhelmed if the number of severe cases becomes very high (refer to Figure~\ref{figure:saturazione} for the present-day situation of intensive care units in Italian regions).

In that context, the rapid filling of the so-called \emph{COVID-19 hospitals} in the areas most affected by the disease represents a serious danger while dealing with sudden emergencies at different scales (villages, cities, regions, countries). Hence, it is compelling to develop new arrangements of the healthcare system, for optimizing the overall structure especially in terms of an efficient distribution of patients among the most suitable facilities.
\smallskip

In this article, we present innovative intervention options, by comparing the standard organisation of the Italian healthcare system with a novel proposal based on the \emph{Abelian Sandpile paradigm}. We provide a systematic approach for an improved planning of the healthcare system by means of a mathematical model for the \emph{self-organisation of critical issues}, which is capable to achieve some specific optimization inside the hospital management. In particular, we aim at reducing the risks of a local/global collapse of hospitals in times of crisis, while improving their functioning in normal times. More specifically, the basic model is grounded in the framework of \emph{cellular automata} and postulates a large network of links between hospitals in a cooperative style of communication. Unfortunately, such connections are currently very limited, which drastically restricts the possibilities of reallocation for the supernumerary patients, and the strong urgency to enhance the hospital network is an important conclusion of the present analysis.

The content of this article is organized as follows.\\
In Section~\ref{section:abelian}, we introduce the Abelian Sandpile paradigm, with a short description of its basic concepts and properties. Section~\ref{section:healthcare} focuses on the adaptation of such paradigm to build up a novel proposal for the organisation of the healthcare system. In Section~\ref{section:examples}, we provide a gallery of examples aiming at illustrating how the proposed model works in the present context of COVID-19 epidemic. More specifically, we analyze two cases relevant to realistic applications: the central and peripheral outbreaks. The article is concluded by Section~\ref{section:beyond} which includes some generalisations of the model and (provisional) conclusions.
\smallskip

\section*{Abbreviations}

\noindent {\tt ASM} -- Abelian Sandpile Model\\
{\tt SOC} -- Self-Organized Criticality\\
{\tt CAM} -- Cellular Automaton Model\\
{\tt SRH} -- Sandpile with Redistribution to the Hub\\
{\tt SID} -- Sandpiles with Internal Dissipation

	
\section{Abelian Sandpile paradigm}
\label{section:abelian}
\smallskip

The notion of Self-Organized Criticality (SOC) was originally introduced by Bak, Tang and Wiesenfeld~\cite{BakTangWies87,BakTangWies88} starting from a basic example proposed as a model for sandpiles (we refer to~\cite{Bak96,Frig03,MarkGros14} for a general introduction to this broad subject).
Since then, the concept has been expanded in many different directions, spanning from classical topics of physics (sandpile avalanches~\cite{Hear03}, distribution of earthquakes~\cite{MalaTurc99,Turc99}, amplitude of solar flares~\cite{MoraChar10}) to less standard economic and socio-political contexts~\cite{AlekHolyKoss02,BartLeinThom05,BionPlucRapi15,ChenZhou08,GlintEtAl10,MaurEtAl18, RamoSassPiqu11,TadiDankMeln17,XiOrmeWang12}, passing through computer networks and biological applications~\cite{Adam95,Newm97,YuanRenShan00,deArPerrHerm06}. At the same time, a huge effort to extend the mathematical tools to deal with theoretical questions has been made, thus contributing to drive SOC into an extraordinary crossroads of probabilistic approaches, graph theory, algebraic geometry, mathematical analysis and optimisation~\cite{Barb17,BasuMoha14,CajuAndr10a,CarlSwin95, CasaZerb00,GaveSchul91,HoffPayt18,PaczBass00}.


\subsection{SOC and Sandpiles}
\label{section:sandpiles}
\smallskip

The original application concerns with modeling of sandpiles, which are regarded as a manageable prototype of SOC, and fundamental contributions have been made by Dhar~\cite{Dhar90,Dhar92,Dhar99} notably in showing the crucial property of \emph{commutativity}. On this account, the adjective Abelian has since been added to the technical terminology, leading to the actual meaning of Abelian Sandpile Model (ASM). Nevertheless, as usual in the most active fields of science, the terminology is not unique, and indeed similar topics are explored under different names (avalanches, chip-firing games, forest-fire models, parking functions, probabilistic abacus, Rotor-Router or Eulerian-walker models, ...) with more or less the same meaning. We refer to~\cite{Creu91,Dhar06,IvasPrie98,Jara18,Mann99} for introductory presentations of the ASM and its various applications.

We are particularly interested to the \emph{load balancing} property, which denotes the method of distributing a set of tasks across a group of resources with the purpose of making their overall processing more efficient by balancing the workload of each operating unit. In the abstract formulation, units are represented by vertexes/nodes of a graph/network, with the corresponding connections represented by edges/links. The objective is to balance the loading process by allowing nodes to exchange \emph{particles} with their neighbors through the incident edges.\\
For the problem covered in this article, we build a specific ASM with \emph{critical height} provided by the capacity of the medical facilities, and we propose to apply the Abelian Sandpile paradigm to achieve a methodical and efficient distribution of patients across the hospitals, in order to maintain the occupancy below the saturation level for dealing with eventual sudden and unexpected emergencies.
\smallskip

A sandpile is also a type of Cellular Automaton Model (CAM), namely a discrete mathematical system fulfilling the following essential conditions:
\begin{itemize}
\item[1.] the evolution takes place without external interventions;
\item[2.] the overall structural development depends only on local rules.
\end{itemize}
To a large extent, the CAM is capable of simulating the dynamics of complex systems which are disposed to organize themselves through unstable critical states until reaching a stable configuration. In addition, the CAM implementation makes it possible to generate global coherent patterns starting from suitable local
instructions, without any external supervisor in charge for understanding the process in its entirety~\cite{Kari05,Sark00,Wolf02}.\\
In practical applications, each restricted portion of space contains a finite number of cells, which assume a finite set of (time-dependent) states. The initial configuration $\bar{\mathbf{z}}^0$ is the combination of a \emph{ground state} configuration $\mathbf{z}^0$ with some perturbation $\mathbf{w}$, and for instance
\begin{equation*}
	\bar{\mathbf{z}}^0 = \mathbf{z}^0 + \mathbf{w}\,.
\end{equation*}
After a given time interval $\Delta t>0$, the system typically comes to a new state $\mathbf{z}^1$, designated \emph{final configuration}, which is determined by the changes of state of the single cell together with all the others, according to the preliminarily established series of (local) rules. Therefore, the evolution is described by means of a function $\Phi$ mapping some input $\bar{\mathbf{z}}^0$ defined (from $\mathbf{z}^0$ and $\mathbf{w}$) over the graph/network $\Gamma$ to an output $\mathbf{z}^1$, as specified by the simple formula $\mathbf{z}^1=\Phi(\bar{\mathbf{z}}^0)$.
\smallskip

It is worth stressing that, in order to properly manage the patients between hospitals of the healthcare system, an optimal underlying ethical-cooperative paradigm has to be stipulated~\cite{Lind17}.



\subsection{Sandpiles on a general network}
\label{section:network}
\smallskip

We start by recollecting some basic definitions in \emph{graph theory}~\cite{Trud94,GrosYell14}.

\begin{definition}
A {\bf graph} (or {\bf network}) is a pair $\Gamma=(X,E)$ where $X$ is a set whose elements are called {\bf vertexes} (or {\bf nodes}) and $E$ is a set of paired vertexes called {\bf edges} (or {\bf links}). A {\bf rooted graph} (or {\bf pointed graph}) is a graph in which one vertex, denoted by $x_\ast$, is distinguished as the {\bf root} of the graph.
\end{definition}

In what follows, the root $x_\ast$ of the graph is also referred to as the {\bf hub}, coherently with the subsequent meaning of the specific application to the healthcare system.

\begin{definition}
An {\bf adjacent vertex} $y\in X$ of $x\in X$ is a vertex which is connected to $x$ by an edge, so that $(x,y)\in E$. The {\bf neighbourhood} $I(x)$ of a vertex $x$ is the subgraph composed of the vertexes adjacent to $x$ and all the edges connecting vertexes adjacent to $x$. The {\bf degree} (or {\bf valency}) of a vertex $x$, denoted by $\textrm{\rm deg}(x)$, is the number of edges which are incident to the vertex $x$.
\end{definition}



\begin{figure}[htb]
\includegraphics[width=0.85\textwidth]{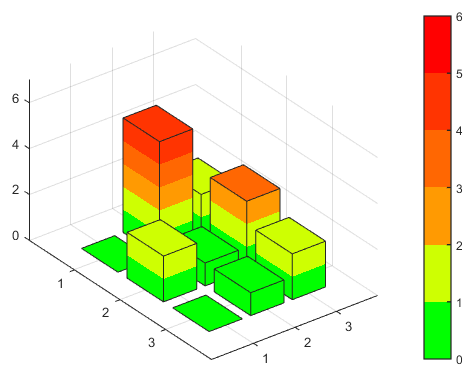}
\caption{\small A simple sandpile with variable height (on a two-dimensional Cartesian grid)}
\label{figure:sandpile}
\end{figure}

Let the graph $\Gamma$ be connected, finite (with a finite number of vertexes denoted by $x_1,x_2,\dots,x_p$), simple (i.e. there are no loops connecting any vertex with itself) and undirected (i.e. the edges are bidirectional). We assume that a stock of identical particles is initially allocated at any vertex $x_i$ of the graph $\Gamma$. The configurations $\mathbf{z}^0$, $\bar{\mathbf{z}}^0$ and $\mathbf{z}^1$, collecting the number of elements located at $x_i$ for any $i$, are natural-valued functions defined on the graph. The rules of the evolution are established to guarantee \emph{stability} of the sandpile dynamics, starting from the definition of the {\bf height function} given by $\mathbf{z} \,:\, X\rightarrow \mathbb{N}^p$ (refer to Figure~\ref{figure:sandpile} for the graphical representation of a sandpile on a two-dimensional Cartesian grid).

\begin{definition}
The height $z_i=\mathbf{z}(x_i)$ is {\bf unstable} (at $x_i\in X$) if $z_i\geq \textrm{\rm deg}(x_i)$, where $\textrm{\rm deg}(x_i)$ is the degree of the vertex $x_i$. Otherwise, the vertex is {\bf stable}.\\
A {\bf stable configuration} is a configuration of vertexes which are stable for any index $i$. A {\bf maximum stable configuration} (or {\bf minimally stable state}) is a configuration in which all the vertexes have a height $z_i$ of one unit lower than the corresponding threshold value $\textrm{\rm deg}(x_i)$. An {\bf almost stable configuration} is a configuration of vertexes which are stable for any index $i$ different from the root, and a {\bf maximum almost stable configuration} is defined accordingly.
\end{definition}

In that framework, an unstable vertex is forced to \emph{topple} over its neighbouring vertexes, thus causing a change of state in the entire configuration.

\begin{definition}
A {\bf toppling} (or {\bf firing}) from the vertex $x_i$ is the mapping $\Xi$ from the graph $\Gamma$ to itself determined by the following (local) rules:
\begin{equation*}
	\Xi(\mathbf{z})_j = z_j - \Delta_{i\to j} \qquad \textrm{with} \quad
	\Delta_{i\to j} := \left\{\begin{aligned}
		& \textrm{\rm deg}(x_i)	&& \textrm{if}\;\;x_j=x_i\\
		& -1					&& \textrm{if}\;\; x_j\in I(x_i)\\
		& \;\; 0				&& \textrm{otherwise}
				\end{aligned}\right.
\end{equation*}
and $\Delta_{i\to j}$ is called the {\bf toppling matrix}~\cite{Dhar90}.
\end{definition}

In principle, different stability criteria could be proposed involving values other than the degree $\textrm{\rm deg}(x_i)$ of the vertex $x_i$, and modified redistribution criteria could also be considered, but we focus on the above definitions for the sake of simplicity in presentation.

The addition of a particle to the sandpile is schematized by summing to the ground state configuration $\mathbf{z}^0=(z_1^0,z_2^0,\dots,z_p^0)$ the vector $\boldsymbol{\delta}_i$ with $1$ at some fixed index $i$ and zeros elsewhere. By iterating this procedure $m$ times for possibly different choices of the index, the initial configuration $\bar{\mathbf{z}}^0$ is finally given by $\mathbf{z}^0+\mathbf{w}$ with $\mathbf{w}=\boldsymbol{\delta}_{i_1}+\boldsymbol{\delta}_{i_2}+\dots+\boldsymbol{\delta}_{i_m}$. It is a straightforward consequence that the operation of adding particles to a sandpile is associative. Moreover, the final configuration $\mathbf{z}^1$ resulting from the evolution is independent of the order in which the topplings are performed, and therefore the operation of toppling is commutative~\cite{Dhar90}.\\
As an illustrative case, we assume that $\textrm{\rm deg}(x_i)=\textrm{\rm deg}(x_j)$ for two adjacent vertexes $x_i$ and $x_j$. Then, we consider a sandpile in which both $x_i$ and $x_j$ are at their critical height, that is $z_i\geq \textrm{\rm deg}(x_i)$ and $z_j\geq \textrm{\rm deg}(x_j)$. According to the dynamics described above, a toppling from the vertex $x_i$ causes the vertex $x_j$ to become unstable, and viceversa, and the sandpile comes to a configuration in which the height $z_k$ for some $k\in\{1,2,\dots,p\}$ increases by a value $\bigtriangleup_{i\to k}+\bigtriangleup_{j\to k}$. Such procedure being symmetrical, by applying this argument repeatedly, we conclude that the same stable configuration is reached irrespective of whether $x_i$ or $x_j$ is toppled first, and in general regardless of the sequence in which unstable vertexes are toppled.


\subsection{Two-dimensional Cartesian grids}
\label{section:cartesian}
\smallskip

\begin{figure}[htb]
\includegraphics[width=0.35\textwidth]{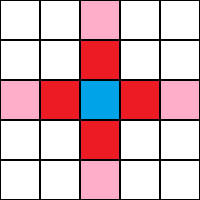}
\qquad\qquad
\includegraphics[width=0.35\textwidth]{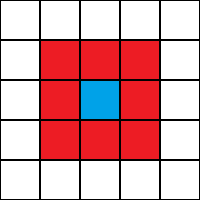}
\caption{\small The von Neumann neighbourhood (left) of the (blue) central cell (an extended neighbourhood includes also the pink cells) and the Moore neighbourhood (right)}
\label{figure:neighbourhood}
\end{figure}

\begin{figure}[htb]
\includegraphics[width=0.85\textwidth]{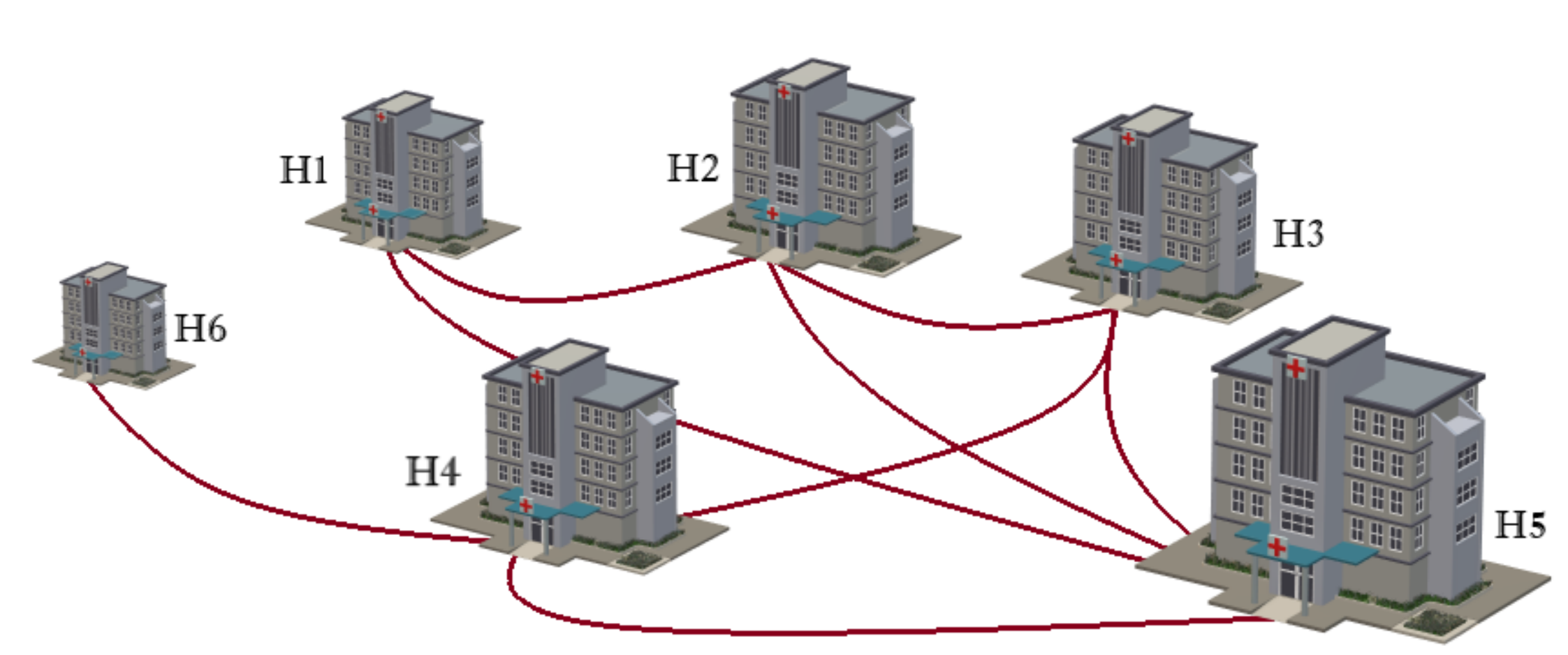}
\caption{\small Arrangement of hospitals on a general unstructured and undirected network.}
\label{figure:network}
\end{figure}

We focus on the simple case of a two-dimensional Cartesian grid, with two standard schemes for adjacent vertexes given by the \emph{von Neumann neighbourhood}, consisting of the four cells obtained moving one step towards North/East/South/West (refer to Figure~\ref{figure:neighbourhood}(left)), and the \emph{Moore neighbourhood}, which includes also the four diagonal cells (refer to Figure~\ref{figure:neighbourhood}(right)). Of course, selecting a more general network provides a higher degree of flexibility and allows to closely represent the connections between various medical facilities of the healthcare system (refer to Figure~\ref{figure:network}).

We consider a special grid of $p=n^2$ cells organized over a squared matrix of size $n\times n$, where the particles are dropped randomly being allowed to stack on top of each other, so that a configuration of the sandpile is described by an element of the space of natural-valued matrices
\begin{equation*}
	\mathcal{M}_{n}(\mathbb{N}) := \bigl\{ A\in \mathbb{R}^{n\times n} \,:\, a_{ij}\in \mathbb{N} \;\; \textrm{for any} \;\; i,j=1,2,\dots,n \bigr\}.
\end{equation*}
The position of each vertex $x_i$ (with its corresponding height $z_{i}\in\mathbb{N}$) is determined by the index $i\in\{1,2,\dots,p=n^2\}$ according to the following common practise:
\begin{equation*}
	\begin{pmatrix}
		1			& 2			& 3		& \dots	& n \\
		n+1			& n+2		& n+3	& \dots	& 2n	\\
		2n+1			& 2n+2		& 2n+3	& \dots	& 3n	\\
		\vdots		& \vdots		& \vdots 	& \ddots	& \vdots \\
		\; n(n-1)+1		& n(n-1)+2	& \dots	& \hdots	& n^2
	\end{pmatrix}
\end{equation*}
For the von Neumann neighborhood, we have $\textrm{\rm deg}(x_i)=4$ for any index $i$ associated with an element of the bulk of the matrix, and for instance
\begin{equation*}
	I(x_{n+2}) = \bigl\{ x_2 \,\textrm{(North)}, \,x_{n+3} \,\textrm{(East)}, \,x_{2n+2} \,\textrm{(South)}, \,x_{n+1} \,\textrm{(West)} \bigr\}.
\end{equation*}
Similarly, for the Moore neighborhood, we have $\textrm{\rm deg}(x_i)=8$ for any index $i$ associated with an element of the bulk of the matrix, and for instance
\begin{equation*}
	\begin{split}
		I(x_{n+2}) = \bigl\{ & x_2 \,\textrm{(North)}, \,x_3 \,\textrm{(North-East)}, \,x_{n+3} \,\textrm{(East)}, \,x_{2n+3} \,\textrm{(South-East)},\\
		& \,x_{2n+2} \,\textrm{(South)}, \,x_{2n+1} \,\textrm{(South-West)}, \,x_{n+1} \,\textrm{(West)}, \,x_1 \,\textrm{(North-West)} \bigr\}.
	\end{split}
\end{equation*}


\begin{example}
\rm We analyze the simple case $n=3$ ($p=n^2=9$) with the von Neumann neighbourhood. By combining the perturbation $\mathbf{w}=4\,\boldsymbol{\delta}_5$ with the ground state configuration $\mathbf{z}^0=\mathbf{0}$, we obtain the following initial and final configurations
\begin{equation*}
	\bar{\mathbf{z}}^0 = \mathbf{0} + 4\,\boldsymbol{\delta}_5
		= \begin{pmatrix} 0 & 0 & 0 \\ 0 & \mathbf{4} & 0 \\ 0 & 0 & 0 \end{pmatrix}
	\quad \longrightarrow \quad
	\mathbf{z}^1 = \Xi(\bar{\mathbf{z}}^0) = \begin{pmatrix} 0 & 1 & 0 \\ 1 & 0 & 1 \\ 0 & 1 & 0 \end{pmatrix}
\end{equation*}
As expected, the threshold value $\textrm{\rm deg}(x_5)=4$ is reached at the central cell, thus inducing a destabilisation of the corresponding vertex, with the consequence that four particles are poured into the adjacent vertexes --the elements of the von Neumann neighbourhood-- to compose a stable final configuration.
\end{example}

As a matter of fact, the toppling from an unstable vertex changes the state of the adjacent vertexes, sometimes causing the appearance of instabilities at some other vertex. Then, subsequent topplings over the adjacent vertexes are generated, triggering a sequence of events which are evocative of \emph{sandpile avalanches} and stopping only when all the cells return strictly below their threshold capacity.\\
If the graph is infinite and connected, with a finite number of particles, all vertexes become stable after a finite number of topplings; moreover, it can be proven that the final (stable) configuration $\mathbf{z}^1$ depends solely on the initial configuration $\bar{\mathbf{z}}^0$, independently of the order in which the topplings are performed~\cite{Dhar90}.\\
On the other hand, if the graph is finite, appropriate boundary conditions must be implemented: assuming that the particles are evacuated from the boundaries, an analogous result of stability holds~\cite{Creu91}, but other types of boundary conditions may not guarantee the same property. During the operation of toppling, no particles are created as they are redistributed to neighbouring cells. For dissipative boundary conditions, the toppling can actually lead to the loss of particles if it occurs on a boundary cell, and we refer to this case as \emph{open boundary conditions}.
\smallskip


\section{Application to healthcare system management}
\label{section:healthcare}
\smallskip

In this article, we pursue the idea of applying the Abelian Sandpile paradigm to the management of patients among the COVID-19 hospitals in Italy.\\
The current organisation assigns responsibility for the healthcare system to the regions, which are restricted territorial bodies with their own statutes, powers and functions established by the Italian Constitution. We refer to the Lazio region when selecting (the order of magnitude of) the number of medical facilities relevant to the numerical simulations. More specifically, the healthcare system of the Lazio region is composed of about $100$ hospitals, mostly located in the metropolitan area of Rome, which are structured within a network of reciprocal connections~\cite{LeoEtAl16}.

In that framework, the {\bf height function} of the nodes (hospitals) indicates the number of hospitalised patients, and the links between different nodes represent the possibility of a direct exchange of patients (that is subject to constraints of geographical proximity, together with inherent organisational requirements of the connecting medical facilities). The node of the network exhibiting the largest number of links is the {\bf hub}, which usually takes on a strategic function for the whole system, whilst the nodes with fewer links identify the hospitals with unfavorable geographical location and reduced capacity/functionality (each hospital having different resources, in the general case the threshold value changes from node to node).


\subsection{Reassignment of outgoing particles}
\label{section:outgoing}
\smallskip

We have already discussed about the crucial role played by the boundary conditions associated with the ASM settled on finite networks, and indeed the existence of stable configurations is strongly influenced by the presence of a \emph{dissipation mechanism} such as the open boundary conditions.\\
Because the loss of particles through the edge of the network and the hypothesis of an infinite network are unrealistic assumptions for practical applications, we have to implement a suitable \emph{reassignment law} for particles generating a critical height at some boundary cells. Hence, we are induced to redistributing the outgoing particles to the root of the network, namely the hub, which is typically a collecting node of the healthcare system. Therefore, we call the resulting process a Sandpile with Redistribution to the Hub (SRH), and we notice that actually the SRH is equivalent to the ASM when there are no topplings occurring at the edge of the network.

In particular, for two-dimensional Cartesian grids with odd size, we incorporate the reassignment of outgoing particles from cells toppling at the border of the matrix by proposing that they are reallocated to the central cell. Furthermore, in case the threshold capacity is reached in several cells at the same time, we impose that the first cell to topple is always the hub, and then additional toppling is performed from the other nodes. In particular, the hub is allowed to topple only once, at the very beginning of the toppling procedure.
\smallskip

We summarize the essential steps of the \emph{CAM workflow} --which turns out to be also the implementation scheme of the numerical codes-- as follows.
\begin{itemize}
\item[1.] {\it Initialisation}.
	We choose a (stable) ground state $\mathbf{z}^0=(z^0_1,z^0_2,\dots,z^0_p)$
	which satisfies the condition that $z^0_i<\textrm{\rm deg}(x_i)$ for any index $i$\,;
\item[2.] {\it Inflow}.
	We combine the ground state $\mathbf{z}^0$ with a perturbation $\mathbf{w}=(w_1,w_2,\dots,w_p)$
	which indicates the number of new patients requiring hospitalization, and we obtain the initial
	configuration $\bar{\mathbf{z}}^0=\mathbf{z}^0+\mathbf{w}$\,;
\item[3.] {\it Hub toppling}.
	If $z^0_\ast\geq \textrm{\rm deg}(x_\ast)$, the hub performs a toppling towards the adjacent nodes
	of its neighbourhood $I(x_\ast)$\,; 
\item[4.] {\it Additional toppling}.
	If $z^0_i\geq \textrm{\rm deg}(x_i)$ for some $x_i\neq x_\ast$, these nodes perform (a sequence
	of) topplings until reaching an almost stable configuration $\mathbf{z}^1=\Phi(\bar{\mathbf{z}}^0)$\,;
\item[5.] {\it Iteration}.
	We reinitialize the ground state with $\mathbf{z}^0$ equal to $\mathbf{z}^1$, and we
	restart from 2.
\end{itemize}


\subsection{Comparison with the standard healthcare system management}
\label{section:comparison}
\smallskip

In order to to compare the efficiency of the novel proposal with the standard organisation of the Italian healthcare system, we provide a reformulation of the current management of connections between hospitals in terms of the CAM paradigm.\\
Perhaps surprisingly, the connection between medical facilities is presently mainly determined by their geographical proximity, corresponding to the basic \emph{adjacency} of nodes. If a patient comes to a hospital where there are no places available, because the threshold capacity has been reached, then the single patient is reallocated to the nearest hospital with available places, and the reassignment is limited to the patients exceeding the threshold value. That being the case, the redistribution of patients is handled manually at each hospitalization, without foreseeing the possibility of repeated similar events which are instead highly probable during sudden and unexpected emergencies.
\smallskip

We summarize the essential steps of the \emph{standard workflow} as follows, by assuming the reassignement law from the edge of the network at the hub as in Section~\ref{section:outgoing}.
\begin{itemize}
\item[1-2.] {\it Initialisation/Inflow}.
	We proceed as in Section~\ref{section:outgoing}\,;
\item[3.] {\it Hub redistribution}.
	If $z^0_\ast\geq \textrm{\rm deg}(x_\ast)$, the hub reallocates only the exceeding patients to the 
	adjacent medical facilities, starting from the less crowded ones (in case of equal crowding, a random choice is operated)\,;
\item[4.] {\it Additional redistribution}.
	If $z^0_i\geq \textrm{\rm deg}(x_i)$ for some $x_i\neq x_\ast$, the redistribution procedure is repeated
	for these nodes, until reaching an almost stable configuration $\mathbf{z}^1=\Psi(\bar{\mathbf{z}}^0)$\,;
\item[5.] {\it Iteration}.
	We reinitialize the ground state with $\mathbf{z}^0$ equal to $\mathbf{z}^1$, and we restart from 2.
\end{itemize}
\smallskip

We illustrate the comparison by considering the following simple example.
\begin{example}
\rm Let the healthcare system be represented by a two-dimensional Cartesian grid of size $3\times 3$, with the von Neumann neighbourhood (so that $\textrm{\rm deg}(x_i)=4$ for any index $i$) and reassignment of outgoing particles to the hub located at the central cell. We choose
\begin{equation*}
	\mathbf{z}^0 = \begin{pmatrix} 2 & 1 & 3 \\ 1 & 3 & 1 \\ 1 & 0 & 2 \end{pmatrix}
	\quad\textrm{and}\quad
	\mathbf{w} = \begin{pmatrix} 0 & 0 & 0 \\ 0 & 1 & 0 \\ 0 & 0 & 0 \end{pmatrix},
\end{equation*}
so that 
\begin{equation*}
	\bar{\mathbf{z}}^0 = \mathbf{z}^0+\mathbf{w}
	= \begin{pmatrix} 2 & 1 & 3 \\ 1 & \mathbf{4} & 1 \\ 1 & 0 & 2 \end{pmatrix}
\end{equation*}
and the initial configuration is unstable. The standard organization manages this critical situation by moving a single particle (patient) from the central cell towards one in the von Neumann neighbourhood, and preferably the South cell, to obtain
\begin{equation*}
	\Psi(\bar{\mathbf{z}}^0) = \begin{pmatrix} 2 & 1 & 3 \\ 1 & 3 & 1 \\ 1 & 1 & 2 \end{pmatrix},
\end{equation*}
where $\Psi$ denotes the \emph{evolution mapping} of the grid according to the standard approach. Then, we reinitialize the ground state with $\Psi(\bar{\mathbf{z}}^0)$ and the perturbation $\mathbf{w}$ is chosen as above, so that the configuration $\Psi(\bar{\mathbf{z}}^0)+\mathbf{w}$ is also unstable with respect to the central cell. By selecting randomly one of the neighbouring cells, another single particle is transferred, and after two additional iterations, the system reaches the final configuration
\begin{equation*}
	\Psi(\bar{\mathbf{z}}^0) = \begin{pmatrix} 2 & 2 & 3 \\ 2 & 3 & 2 \\ 1 & 1 & 2 \end{pmatrix}.
\end{equation*}
On the other hand, the ASP paradigm suggests to transfer at the same time all the four patients initially allocated to the hub towards the von Neumann neighbourhood, to obtain
\begin{equation*}
	\Phi(\mathbf{z}^0) = \begin{pmatrix} 2 & 2 & 3 \\ 2 & 0 & 2 \\ 1 & 1 & 2 \end{pmatrix}.
\end{equation*}
Finally, under the assumption that new patients are always collected at the hub, the system can accept exactly $3$ more particles from the additional iterations, before reaching a critical situation as before. Of course, it could be argued that, being the final configurations the same, there would be no obvious evidence to prefer the novel approach to the standard organization. However, the deciding factor is that we have unified within a single time-step the efforts for reallocating patients among the medical facilities, which is usually regarded as a source of stress for the overall healthcare system. Hence, we have left the structure time to reorganize without suffering a constant state of emergency, which is precisely the issue of the so-called \emph{predictive logistics}. 
\end{example}


\section{Examples and numerical simulations}
\label{section:examples}
\smallskip

We collect numerical results from a selection of preliminary cases, which are nevertheless useful to understand the alternative methods of healthcare system management provided by the two approaches described in Section~\ref{section:healthcare}, and their inherent dynamics. We consider two-dimensional Cartesian grids, which are arranged into $n\times n$ (squared) matrices with the Moore neighbourhood, so that $\textrm{\rm deg}(x_i)=8$ for any index $i\in\{1,2,\dots,p=n^2\}$, and we assume the reassignment to the hub of outgoing particles from the edge of the network as in Section~\ref{section:outgoing}.\\
We recall that each cell/node of the grid represents a hospital, and the height function reproduces the number of patients in the medical facilities, expressed as the percentage of capacity already achieved.

We start by focusing on the possible outputs of a basic case consisting in a single hospital originally attaining its threshold capacity, which is located at the hub (central cell), with a (randomly chosen) stable ground state configuration. To highlight possibly critical situations, we boldmark all values above the threshold capacity, and we encircle the values $5$, $6$ and $7$ dangerously close to the saturation level.


\begin{example}
\rm For $n=3$, the hospitals are $p=n^2=9$, and we choose an initial configuration with all new patients coming to the hub $x_{5}$, given by the perturbation $\mathbf{w}=4\,\boldsymbol{\delta}_5$, so that
\begin{equation*}
	\mathbf{z}^0 = \begin{pmatrix} 2 & 3	& 1 \\ \circled{5} & \circled{7} & 2 \\ 4 & 3 & 3 \end{pmatrix}
	\quad\textrm{and}\quad
	\bar{\mathbf{z}}^0 = \begin{pmatrix} 2 & 3 & 1 \\ \circled{5} & {\bf 11} & 2 \\ 4 & 3 & 3 \end{pmatrix}.
\end{equation*}
Because the hub has to manage a number of hospitalisations greater than its threshold capacity, some patients must be transferred to adjacent facilities.\\
According to the standard strategy, only four patients need a different allocation to be found among the neighbourhing cells, ending at the (stable) configuration
\begin{equation*}
	\Psi(\bar{\mathbf{z}}^0)
	= \begin{pmatrix}	2		& 4 		& 1	\\
					\circled{5}		& \circled{7} 	& 4	\\
					4 		& 4		& 3	\end{pmatrix}.
\end{equation*}
As it is clearly seen, an additional iteration with $\mathbf{w}=4\,\boldsymbol{\delta}_5$ determines a new unstable/critical configuration at the central node (with again four patients in excess), demanding for a further reorganisation of the hub by means of another transfer operation. On the other hand, by employing the SRH approach, we obtain a different outcome given by
\begin{equation*}
	\Phi(\bar{\mathbf{z}}^0)
	= \begin{pmatrix}	3	 & 4	& 2	\\
					\circled{6} & 3	& 3	\\
					\circled{5}	 & 4	& 4	\end{pmatrix},
\end{equation*}
such situation being more favorable in terms of load-balancing since an additional iteration with perturbation $\mathbf{w}=4\,\boldsymbol{\delta}_\ast$ does not produce any unstable configuration.
Moreover, there is still a certain amount of available places at the hub, and therefore the novel proposal solves the criticalities possibly generated by an emergency.
\end{example}

Despite its simplicity, the above example suggests to attribute to each configuration the value of a given functional, aiming to provide an easily readable \emph{indicator} of the effectiveness of the allocation strategy. This is necessarily a very delicate issue.\\
As a first attempt, a possible choice is to count the total number of medical facilities attaining a given fraction of their threshold capacity, which are denoted by \emph{critical points}. Such choice is actually quite questionable, since it does not take keep memory of the value of incoming patients at the
beginning of the iteration.\\
A different (rough) quantitative measure of the system management efficiency is
\begin{equation*}
	\mathcal{F}(\mathbf{w},\mathbf{z}) := \frac{\mathbf{w}}{\sum_{j} w_j}\cdot \mathbf{z}\,,
\end{equation*}
which intends to evaluate the risk that a new patient comes to a given location, taking into account the previous history of the system (in this case, the effect of the perturbation $\mathbf{w}$).
\smallskip

For the case illustrated in the Example~3, the function $\mathcal{F}$ equals the state of the hub $x_5$, 
so that
\begin{equation*}
	\mathcal{F}(\mathbf{w},\bar{\mathbf{z}}^0) = \bar z^0_5 = 11,\quad
	\mathcal{F}\bigl(\mathbf{w},\Psi(\bar{\mathbf{z}}^0)\bigr) = \Psi(\bar{\mathbf{z}}^0)_5 = 7,\quad
	\mathcal{F}\bigl(\mathbf{w},\Phi(\bar{\mathbf{z}}^0)\bigr) = \Phi(\bar{\mathbf{z}}^0)_\ast = 3\,.
\end{equation*}
The minimum value is achieved for the configuration $\Phi(\bar{\mathbf{z}}_0)$ corresponding to the SRH approach (or, equivalently, to ASM).


\begin{example}
\rm For $n=5$, the hospitals are $p=n^2=25$, and we choose an initial configuration with all new patients coming to the hub $x_{13}$, given by the perturbation $\mathbf{w}=\boldsymbol{\delta}_{9}+\boldsymbol{\delta}_{12}+4\boldsymbol{\delta}_13
+2\boldsymbol{\delta}_{14}+\boldsymbol{\delta}_{18}+\boldsymbol{\delta}_{19}$, so that
\begin{equation*}
	\mathbf{z}^0 = \begin{pmatrix}
		1 & 2 & 4		& 2		& \circled{5} \\
		2 & 4 & 2		& 3		& 1 \\
		3 & 2 & \circled{7}	& 2		& 3 \\
		2 & 1 & 4		& 2 		& 2 \\
		4 & 2 & 1		& \circled{5}	& 4
			\end{pmatrix}
	\quad\textrm{and}\quad
	\bar{\mathbf{z}}^0 = \begin{pmatrix}
		1 & 2 & 4 		& 2		& \circled{5} \\
		2 & 4 & 2		& 4		& 1 \\
		3 & 3 & {\bf 11}	& 4		& 3 \\
		2 & 1 & \circled{5}	& 3		& 2 \\
		4 & 2 & 1		& \circled{5}	& 4
			\end{pmatrix}.
\end{equation*}
According to the standard strategy, an admissible solution could be
\begin{equation*}
	\Psi(\bar{\mathbf{z}}^0) = \begin{pmatrix}
		1 & 2 & 4		& 2		& \circled{5} \\
		2 & 4 & 4		& 4		& 1 \\
		3 & 3 & \circled{7}	& 4		& 3 \\
		2 & 3 & \circled{5}	& 3		& 2 \\
		4 & 2 & 1		& \circled{5}	& 4
			\end{pmatrix},
\end{equation*}
and the same observation previously done holds, that is such configuration likely determines an instability located at the hub after a subsequent process iteration, since the number of hospitalised patients is very close to the critical threshold.\\
On the other hand, by employing the SRH approach, we obtain
\begin{equation*}
	\Phi(\bar{\mathbf{z}}^0) = \begin{pmatrix}
		1 & 2 	& 4 		& 2 		&  \circled{5} \\
		2 & \circled{5}	& 3		& \circled{5} 	& 1 \\
		3 & 4		& 3		& \circled{5} 	& 3 \\
		2 & 2		& \circled{6}	& 4 		& 2 \\
		4 & 2		& 1		& \circled{5}	& 4
			\end{pmatrix},
\end{equation*}
where the number of patients hospitalised at the hub is far away from the critical level. Indeed, by computing the indicator function $\mathcal{F}$ at the different outcomes, we deduce that
\begin{equation*}
	\mathcal{F}\bigl(\mathbf{w},\bar{\mathbf{z}}^0\bigr) = 6.7,\qquad
	\mathcal{F}\bigl(\mathbf{w},\Psi(\bar{\mathbf{z}}^0)\bigr) = 5.1,\qquad
	\mathcal{F}\bigl(\mathbf{w},\Phi(\bar{\mathbf{z}}^0)\bigr) = 4.1\,,
\end{equation*}
since $\sum_{j} w_j=10$, and the smallest value is achieved by the configuration resulting from the SRH approach; thus, such strategy must be preferred to the standard one.
\end{example}


Next, we explore a case where two subsequent topplings occur.

\begin{example}
\label{ex:10}
\rm For $n=5$, we choose the ground state and the initial configuration as follows,
\begin{equation*}
	\mathbf{z}^0 = \begin{pmatrix}
		4		& 1	& 0		& 1		& 3		\\
		\circled{5}		& 0		& \circled{5}	& 1	& 1		\\
		1 		& 2		& \circled{7}	& \circled{7}	& 4		\\
		\circled{5}		& \circled{5}	& 2		& 4 		& \circled{5}	\\
		3		& \circled{5}	& 4		& \circled{5}	& 3
			\end{pmatrix}
	\quad\textrm{and}\quad
	\bar{\mathbf{z}}^0 = \begin{pmatrix}
		4		& 1		& 0			& 1		& 3		\\
		\circled{5}		& 0		& \circled{5}		& 1		& 1		\\
		1 		& 2		& {\bf 11} & \circled{7}	& 4		\\
		\circled{5}		& \circled{5}	& 2			& 4 		& \circled{5}	\\
		3		& \circled{5}	& 4			& \circled{5}	& 3
			\end{pmatrix},
\end{equation*}
where $\mathbf{w}=4\boldsymbol{\delta}_{13}$. Then, according to the standard strategy, a final configuration is given by
\begin{equation*}
	\Psi(\bar{\mathbf{z}}^0) = \begin{pmatrix}
		4		& 1		& 0		& 1		& 3		\\
		\circled{5}		& 3		& \circled{5}	& 2		& 1		\\
		1 		& 2		& \circled{7}	& \circled{7}	& 4		\\
		\circled{5}		& \circled{5}	& 2		& 4 		& \circled{5}	\\
		3		& \circled{5}	& 4		& \circled{5}	& 3
			\end{pmatrix}.
\end{equation*}
For applying the SRH approach, we pass through an intermediate configuration
\begin{equation*}
	\begin{pmatrix}
		4		& 1		& 0		& 1			& 3		\\
		\circled{5}		& 0+1	& \circled{5+1}	& 1+1		& 1		\\
		1 		& 2+1	& 11-8	& \mathbf{7+1}	& 4		\\
		\circled{5}		& \circled{5+1}	& 2+1	& \circled{4+1} 	& \circled{5}	\\
		3		& \circled{5}	& 4		& \circled{5}		& 3
			\end{pmatrix}
	=
	\begin{pmatrix}
		4		& 1		& 0		& 1			& 3		\\
		\circled{5}		& 1		& \circled{6}	& 2			& 1		\\
		1 		& 3		& 3		& \mathbf{8}	& 4		\\
		\circled{5}		& \circled{6}	& 3		& \circled{5}	 	& \circled{5}	\\
		3		& \circled{5}	& 4		& \circled{5}		& 3
			\end{pmatrix},
\end{equation*}
in which the cell $x_{12}$ is unstable, and finally we arrive at the final configuration
\begin{equation*}
	\Phi(\bar{\mathbf{z}}^0) = \begin{pmatrix}
		4		& 1		& 0		& 1		& 3		\\
		\circled{5}		& 1		& \circled{6+1}	& 2+1	& 1+1	\\
		1 		& 3		& 3+1	& 8-8	& \circled{4+1}	\\
		\circled{5}	& \circled{6}	& 3+1	& \circled{5+1} & \circled{5+1}	\\
		3		& \circled{5}	& 4		& \circled{5}	& 3
			\end{pmatrix}
	= \begin{pmatrix}
		4		& 1		& 0		& 1		& 3		\\
		\circled{5}	& 1		& \circled{7}	& 3		& 2		\\
		1 		& 3		& 4		& 0		& \circled{5}	\\
		\circled{5}	& \circled{6}	& 4		& \circled{6}	& \circled{6}	\\
		3		& \circled{5}	& 4		& \circled{5}	& 3
			\end{pmatrix}.
\end{equation*}

The matrices $\Psi(\bar{\mathbf{z}}^0)$ and $\Phi(\bar{\mathbf{z}}^0)$ are represented in Figure~\ref{figure:figure6} with different colors assigned to cells, corresponding to the degree of saturation achieved.

\begin{figure}[htb]
\includegraphics[width=0.90\textwidth]{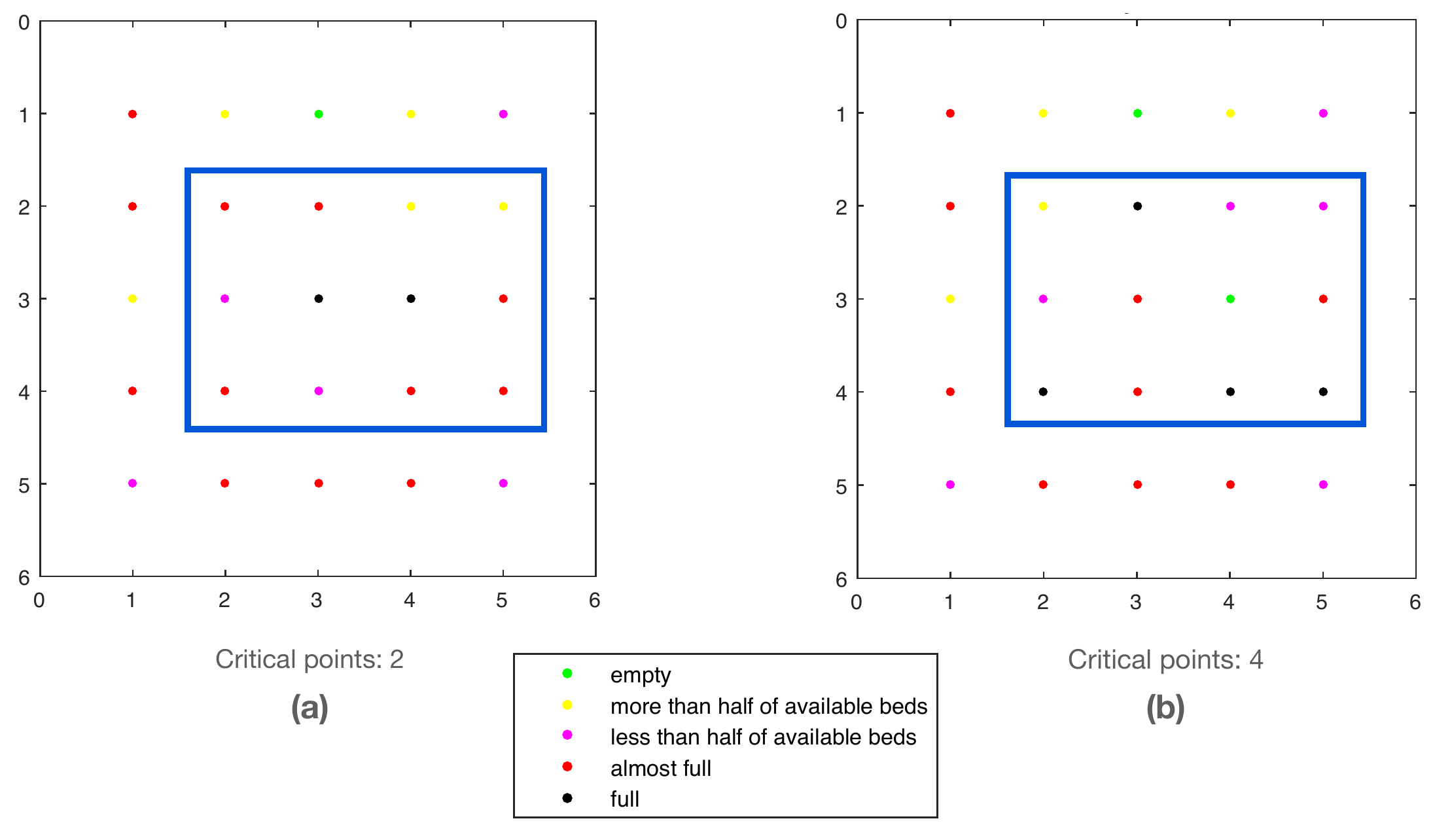}
\caption{\small Model comparison from the Example~\ref{ex:10}: final configurations $\Psi(\bar{\mathbf{z}}^0)$ (left) and $\Phi(\bar{\mathbf{z}}^0)$ (right). The colors attached to each cell/node represent its corresponding relative capacity (0=green, 1/2=yellow, 3/4=magenta, 4/5=red, 6/7=black)}
\label{figure:figure6}
\end{figure}

At the bottom of each subfigure, we report the counting of critical points, given by the number of cells with values $6$ or $7$. Although a quick glance can make us suspect that the left configuration is preferable, the situation is indeed different. In fact, first of all we distinguish the presence of two almost critical values $7$ inside the configuration plotted on the left, against the presence of only one value $7$ from the second approach. Furthermore, the hub --the facility most at risk taking into account recent history-- has a value $7$ in the first case, whilst a value $4$ occurs in the second case, thus making the latter situation preferable.

Finally, we compute the function $\mathcal{F}$ corresponding to the different configurations: since $\mathbf{w}/\sum_j w_j$ is everywhere zero except at the hub, where it is equal to $1$, the values of $\mathcal{F}$ coincide with the values at the hub, which are given by
\begin{equation*}
	\mathcal{F}(\mathbf{w},\bar{\mathbf{z}}^0) = 11,\qquad
	\mathcal{F}\bigl(\mathbf{w},\Psi(\bar{\mathbf{z}}^0)\bigr) = 7,\qquad
	\mathcal{F}\bigl(\mathbf{w},\Phi(\bar{\mathbf{z}}^0)\bigr) = 4,
\end{equation*}
suggesting again the advantages of configurations proposed by the SRH approach.
\end{example}


\subsection{Simulation of multiple central outbreaks}
\label{section:outbreak1}
\smallskip

We analyze the case with more than one hospital in the same restricted area (including the hub) which is concerned with new incoming patients. One might think of this event as the emergence of multiple outbreaks within some smaller district of a larger area. We consider a larger network, with $n=9$, whose size is comparable to the capacity of the healthcare system in the Lazio region in Italy~\cite{LeoEtAl16}.

We consider the perturbation $\mathbf{w}=2\boldsymbol{\delta}_{32}+\boldsymbol{\delta}_{33}+5\boldsymbol{\delta}_{41}+2\boldsymbol{\delta}_{42}$ and the ground state
\begin{equation*}
	\mathbf{z}^0 = \begin{pmatrix}
		1 	 & 2 		& 1 		& 4 		&  \circled{6}	& 2 		& 3 	&  \circled{6}	& 2	\\
		3 	 & 2 		&  \circled{5}	& 2 		& 1	 	& 3 		& 2 		& 4	  	& 3	\\
		3 	 & 3 		& 1 	& \circled{5} 	&  \circled{6}	& 2 		& 3 		& 1		& 3	\\
		 \circled{6} & 1 	&  \circled{5} 	& 3 		& 2 		&  \circled{5} 	& 2 	& 3		& 4	\\
		1 	 & 3 		& 2 	& 1 	&  \circled{6}	& 3 		& 4 	&  \circled{5}	&  \circled{5}	\\
		1 	 & 4 		& 1 		& 2 		& 2 		&  \circled{5}	& 1 		& 2		& 3 	\\
		4 	 & 1 		& 3 		& 4 		&  \circled{5} 	&  \circled{6}	& 2 	&  \circled{5} & 3 \\
		4 	 & 2 		& 3 		&  \circled{5}	& 2 		& 2 		&  \circled{6}	& 3	& 1	\\
		1 	 &  \circled{6}	&  \circled{5} 	& 2 		& 4 		& 4 		& 2 		& 1	& 2
	\end{pmatrix},
\end{equation*}
so that $\sum_{j} z^0_j=250$ and $\sum_{j} w_j=10$, corresponding to the initial configuration
\begin{equation*}
\begin{aligned}
	\bar{\mathbf{z}}^0 &= \begin{pmatrix}
		1 	&2 		& 1	&4		& \circled{6} 	&2		&3		&\circled{6}	&2	\\
		3	&2 		& \circled{5} &2		& 1 		&3		&2		&4		&3		\\
		3	&3 		& 1	& \circled{5}	&  \circled{6} 	&2		&3		&1		&3	\\
		 \circled{6} &1 	&  \circled{5} &3 & 2+2 	& \circled{5+1}	&2	&3		&4	\\
		1	&3 		& 2	&1	& \circled{6+5} & \circled{3+2}	&4	& \circled{5}	& \circled{6} \\
		1	&4 		& 1	&2		& 2 		& \circled{5}	&1 		&2		&3		\\
		4	&1 		& 3	&4		&  \circled{5}	& \circled{6} 	&2 	& \circled{5}	&3	\\
		4	&2 		& 3	& \circled{5}	& 2		&2 		& \circled{6}	&3		&1	\\
		1	& \circled{6}	&  \circled{5} &2	& 4		&4 		&2		&1	&2
	\end{pmatrix}\\ \\
		&= \begin{pmatrix}
		1	&2 &1	&4		& \circled{6}	&2		&3		& \circled{6}	&2	\\
		3	&2 & \circled{5}	&2		&1		&3		&2		&4		& \circled{6} \\
		3	&3 &1	& \circled{5}	& \circled{6}	&2		&3		&1		&3	\\
		 \circled{6} &1 & \circled{5} &3		&4		& \circled{6}	&2		&3	&4	\\
		1	&3 &2	&1		&{\bf 11}	& \circled{5}	&4	& \circled{5}	& \circled{6} \\
		1	&4 &1	&2		&2		& \circled{5}	&1		&2		&3		\\
		4	&1 &3	&4		& \circled{5}	& \circled{6}	&2		& \circled{5}	&3	\\
		4	&2 &3	& \circled{5}	&2		&2		& \circled{6}	&3		&1	\\
		1	&3 & \circled{5} &2		&4		&4		&2		&1		&2
	\end{pmatrix}.
\end{aligned}
\end{equation*}
Then, according to the standard strategy, a final configuration is given by
\begin{equation*}
	\Psi(\bar{\mathbf{z}}^0) = \begin{pmatrix}
		1	& 2 		& 1		& 4     	& \circled{6}    	&2      	&3	& \circled{6}	& 2	\\
		3	& 2 		& \circled{5}	& 2     	&1      	&3      	&2		& 4 		& 3	\\
		3	& 3 		& 1 		& \circled{5}	& \circled{6}     	&2      	&3	& 1 	& 3	\\
		\circled{6} & 1 		& \circled{5} 	& 3+0	& 4+0  	&\circled{6+0}	&2	& 3 	& 4	\\
		1	& 3 	& 2 	& 1+2	& \circled{11-4} &\circled{5+0}	&4	& \circled{5} & \circled{6}	\\
		1	& 4 		& 1 		& 2+1	& 2+1  	&\circled{5+0}	&1		& 2 	& 3	\\
		4	& 1 		& 3 		& 4    	& \circled{5}	&6     	&2	& \circled{5} & 3 \\
		4	& 2 		& 3 		& \circled{5}   	&2       	&2    		&\circled{6}	& 3 	& 1	\\
		1	& \circled{6} 	& \circled{5}	& 2    &4      	&4     	&2 		& 1 		& 2
	\end{pmatrix}.
\end{equation*}
Hence, the central subgraph around the hub is composed by the elements
\begin{equation*}
	\begin{pmatrix} 3 &  4 & \circled{6} \\ 3 & \circled{7} & \circled{5} \\ 3 & 3 & \circled{5} \end{pmatrix}.
\end{equation*}
In the subsequent steps, any additional patients coming to the hub destabilises the configuration.
On the other hand, by employing the SRH approach, we obtain
\begin{equation*}
	\Phi(\bar{\mathbf{z}}^0)=\begin{pmatrix}
		1 	& 2 		& 1		& 4 		& \circled{6}	& 2 		& 3	& \circled{6}	& 2	\\
		3 	& 2 		& \circled{5}	& 2 		& 1 		& 3 		& 2 		& 4		& 3	\\
		3 	& 3 		& 1		& \circled{5} 	& \circled{6}	& 2 	& 3 		& 1		& 3	\\
		\circled{6} & 1 		& \circled{5}	& 3+1 & \circled{4+1} & \circled{6+1}	& 2 	& 3	& 4	\\
		1 	& 3 		& 2	& 1+1	& 11-8	& \circled{5+1}	& 4 	& \circled{5}	& \circled{6} \\
		1 	& 4 		& 1		& 2+1 	& 2+1 	& \circled{5+1}	& 1 	& 2		& 3	\\
		4 	& 1 		& 3		& 4 		& \circled{5} 	& \circled{6}	& 2 	&\circled{5} & 3	\\
		4 	& 2 		& 3 		& \circled{5}	& 2 		& 2		& \circled{6}	& 3	& 1	\\
		1 	& \circled{6}	& \circled{5}	& 2 		& 4 		& 4		& 2 		& 1	& 2
	\end{pmatrix},
\end{equation*}
and the corresponding central subgraph around the hub is
\begin{equation*}
	\begin{pmatrix}	4 & \circled{5} 	& \circled{7} \\ 2 & 3 	& \circled{6} \\ 3 & 3 	& \circled{6}  \end{pmatrix}.
\end{equation*}
The difference between $\Psi(\bar{\mathbf{z}}^0)$  and $\Phi(\bar{\mathbf{z}}^0)$ is transparent regarding,
specifically, the number of patients hospitalised in the hub. A further iteration with the same perturbation $\mathbf{w}$ would lead to a new critical situation for the hub in the first case, but it does not in the second.

\begin{figure}[htb]
\includegraphics[width=0.90\textwidth]{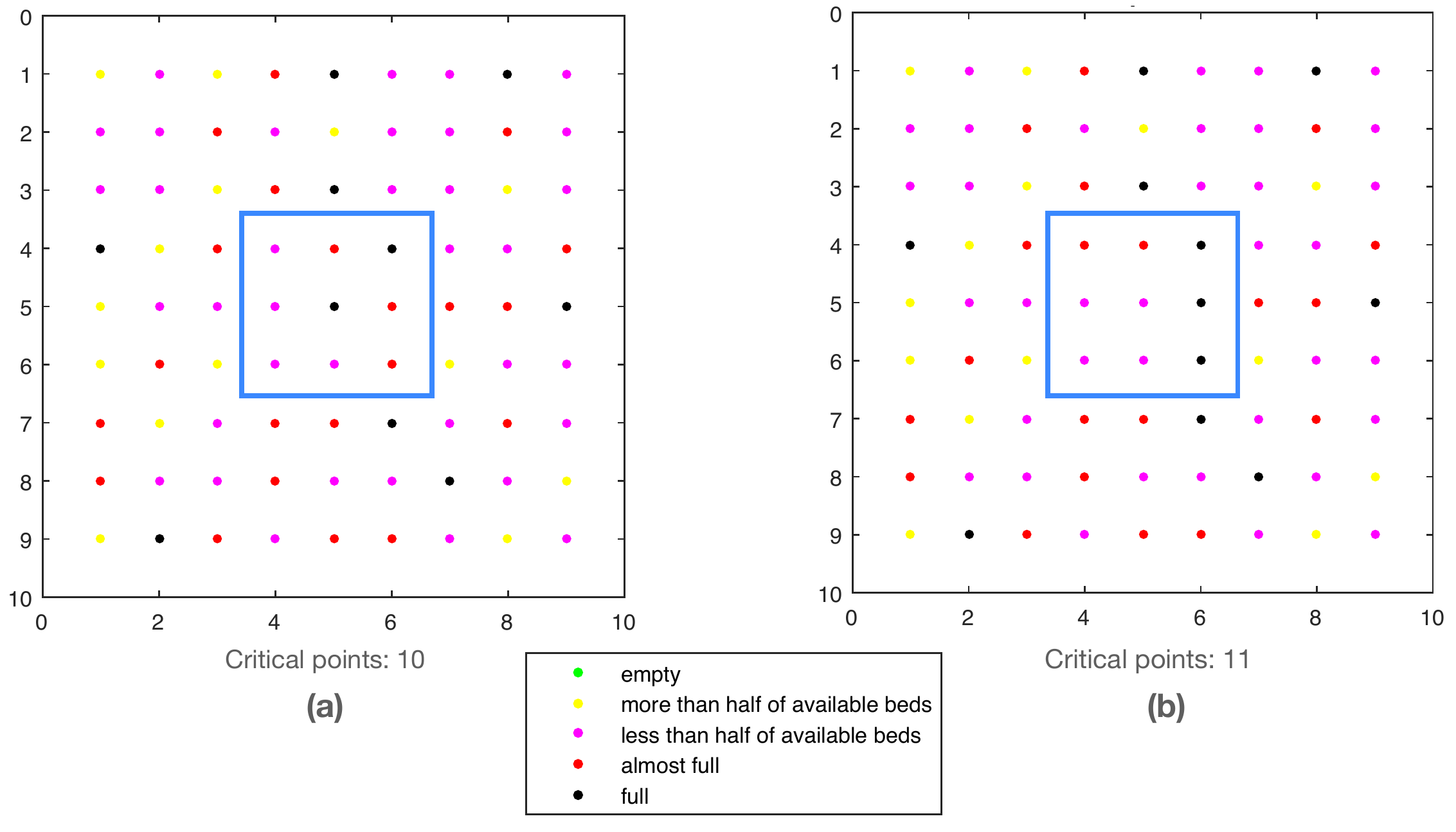}
\caption{\small Model comparison from the case of central areas outbreaks: final configurations $\Psi(\bar{\mathbf{z}}^0)$ (left) and $\Phi(\bar{\mathbf{z}}^0)$ (right). The colors attached to each cell/node represent its corresponding relative capacity (0=green, 1/2=yellow, 3/4=magenta, 4/5=red, 6/7=black)}
\label{figure:figure7}
\end{figure}

Figure~\ref{figure:figure7} provides a representation of the two configurations $\Psi(\bar{\mathbf{z}}^0)$ and $\Phi(\bar{\mathbf{z}}^0\bigr)$.

Counting the critical points in Figure~\ref{figure:figure7} could give the (incorrect) impression that the standard approach has to be preferred with respect to the novel proposal. But considering the indicator function $\mathcal{F}$ applied to $\bar{\mathbf{z}}^0$, $\Psi(\bar{\mathbf{z}}^0)$
and $\Phi(\bar{\mathbf{z}}^0)$, we obtain
\begin{equation*}
	\mathcal{F}\bigl(\mathbf{w},\bar{\mathbf{z}}^0\bigr) =7.9,\qquad
	\mathcal{F}\bigl(\mathbf{w},\Psi(\bar{\mathbf{z}}^0)\bigr)=5.9,\qquad
	\mathcal{F}\bigl(\mathbf{w},\Phi(\bar{\mathbf{z}}^0\bigr))=4.5,
\end{equation*}
as a consequence of the fact that memory of the dynamics --represented by $\mathbf{w}$-- is now considered, and again the minimum value is achieved for the SRH strategy.


\subsection{Simulation of peripheral outbreaks}
\label{section:outbreak2}
\smallskip

Next, we focus on the case of outbreaks occurring in the vicinity of a node located at the edge of the Cartesian grid. Specifically, we take $n=9$ and the same ground state $\mathbf{z}^0$ of the previous example, and an \emph{inflow matrix} of patients $\mathbf{w}$ given by the submatrix
\begin{equation*}
	\begin{pmatrix}  \circled{5} & 2 \\ 2 & 1 \end{pmatrix}
\end{equation*}
located at the top-right corner of the null matrix.
Then, the initial configuration $\bar{\mathbf{z}}^0=\mathbf{z}^0+\mathbf{w}$ is given by
\begin{equation*}
	\bar{\mathbf{z}}^0 = \begin{pmatrix}
		1	& 2	& 1		& 4		&  \circled{6}	& 2 		& 3		& \mathbf{11}	& 4 	\\
		3 	& 2 	&  \circled{5}	& 2 		& 1 	& 3 		& 2		&  \circled{6} 		& 4 	\\
		3 	& 3 	& 1 	&  \circled{5} 	&  \circled{6}	& 2 		& 3		& 1 			& 3 	\\
		 \circled{6} & 1 	&  \circled{5} 	& 3 	& 2 		&  \circled{5} 	& 2		& 3 		& 4	\\
		1 	& 3 	& 2 		& 1 	&  \circled{6}	& 3 		& 4	&  \circled{5} 	&  \circled{6}	\\
		1 	& 4 	& 1 		& 2 		& 2 	&  \circled{5}	& 1		& 2 			& 3	\\
		4 	& 1 	& 3 		& 4 	&  \circled{5}	&  \circled{6} 	& 2 		&  \circled{5} 	& 3	\\
		4 	& 2 	& 3 		&  \circled{5} 	& 2	& 2 	&  \circled{6}	& 3 			& 1	\\
		1 	&  \circled{6} 	&  \circled{5} 	& 2 		& 4		& 4 		& 2 		& 1 	& 2
	\end{pmatrix},
\end{equation*}
with the cell above the threshold capacity located at position $(1,8)$. Therefore, four patients are in excess and require a different reallocation. Among others, a possible output of the standard strategy is the final configuration given by
\begin{equation*}
	\Psi(\bar{\mathbf{z}}^0) = \begin{pmatrix}
		1 	& 2		& 1		& 4	&  \circled{6}	& 2	&  \circled{5}	&  \circled{7} 	& 4	\\
		3 		& 2		&  \circled{5}	& 2 		& 1 		& 3	& 4	&  \circled{6} 	& 4	\\
		3		& 3 		& 1 		&  \circled{5}	&  \circled{6}	& 2	& 3		& 1	& 3	\\
		 \circled{6}		& 1  &  \circled{5}	& 3 	& 2 		&  \circled{5} 	& 2		& 3 	& 4	\\
		1 	& 3 		& 2		& 1 	&  \circled{6}	& 3 		& 4 	&  \circled{5} 	&  \circled{6} \\
		1 	& 4 		& 1		& 2 		& 2 		&  \circled{5} 	& 1 		& 2 		& 3	\\
		4 	& 1 		& 3		& 4 	&  \circled{5} 	&  \circled{6} 	& 2 	&  \circled{5} 	& 3	\\
		4 		& 2 	& 3		&  \circled{5} 	& 2 		& 2 	&  \circled{6}	& 3 	& 1	\\
		1 	&  \circled{6} 	&  \circled{5}	& 2 		& 4		& 4 		& 2	& 1 	& 2
	\end{pmatrix}.
\end{equation*}
On the other hand, by employing the SRH approach, we pass through an intermediate state
\begin{equation*}
	\begin{pmatrix}
		1	& 2	& 1	& 4		&  \circled{6}	& 2 		& 4		& 3	&  \circled{5} 	\\
		3 	& 2 	&  \circled{5}	& 2 	& 1 		& 3 	&  \circled{3}	&  \circled{7} 	&  \circled{5} \\
		3 	& 3 	& 1 	&  \circled{5} 	&  \circled{6}	& 2 		& 3		& 1 	& 3 	\\
		 \circled{6} & 1 	&  \circled{5} 	& 3 		& 2 			&  \circled{5} 	& 2	& 3 	& 4	\\
		1 	& 3 	& 2 		& 1 	& \mathbf{9}	& 3 		& 4		&  \circled{5} 	&  \circled{6} \\
		1 	& 4 	& 1 		& 2 	& 2 		&  \circled{5}	& 1		& 2 	& 3	\\
		4 	& 1 	& 3 		& 4 	&  \circled{5}	&  \circled{6} 	& 2 	&  \circled{5} 	& 3	\\
		4 	& 2 	& 3 		&  \circled{5} 	& 2		& 2 		&  \circled{6}	& 3 	& 1	\\
		1 	&  \circled{6} 	&  \circled{5} 	& 2 		& 4		& 4 	& 2 		& 1 	& 2
	\end{pmatrix},
\end{equation*}
to reach the final configuration
\begin{equation*}
	\Phi(\bar{\mathbf{z}}^0)
		=\begin{pmatrix}
		1	& 2	& 1	& 4		&  \circled{6}	& 2 		& 4	& 3		&  \circled{5} 	\\
		3 	& 2 	&  \circled{5}	& 2 	& 1 		& 3 	&  \circled{3}	&  \circled{7} 	&  \circled{5}  \\
		3 	& 3 	& 1 	&  \circled{5} 	&  \circled{6}		& 2 		& 3		& 1 		& 3 	\\
		 \circled{6} & 1 		&  \circled{5} 	& 4 		& 3 	&  \circled{6} 	& 2	& 3 	& 4	\\
		1 	& 3 	& 2 		& 2 		& 1		& 4		& 4		&  \circled{5} 	&  \circled{6}	\\
		1 	& 4 	& 1 		& 3 		& 3 	&  \circled{6}	& 1		& 2 		& 3	\\
		4 	& 1 	& 3 		& 4 		&  \circled{5}	&  \circled{6} 	& 2 	&  \circled{5} 	& 3	\\
		4 	& 2 	& 3 		&  \circled{5} 	& 2			& 2 		&  \circled{6}	& 3 	& 1	\\
		1 	&  \circled{6} 	&  \circled{5} 	& 2 		& 4	& 4 		& 2 		& 1 	& 2
	\end{pmatrix}.
\end{equation*}

\begin{figure}[htb]
\includegraphics[width=0.90\textwidth]{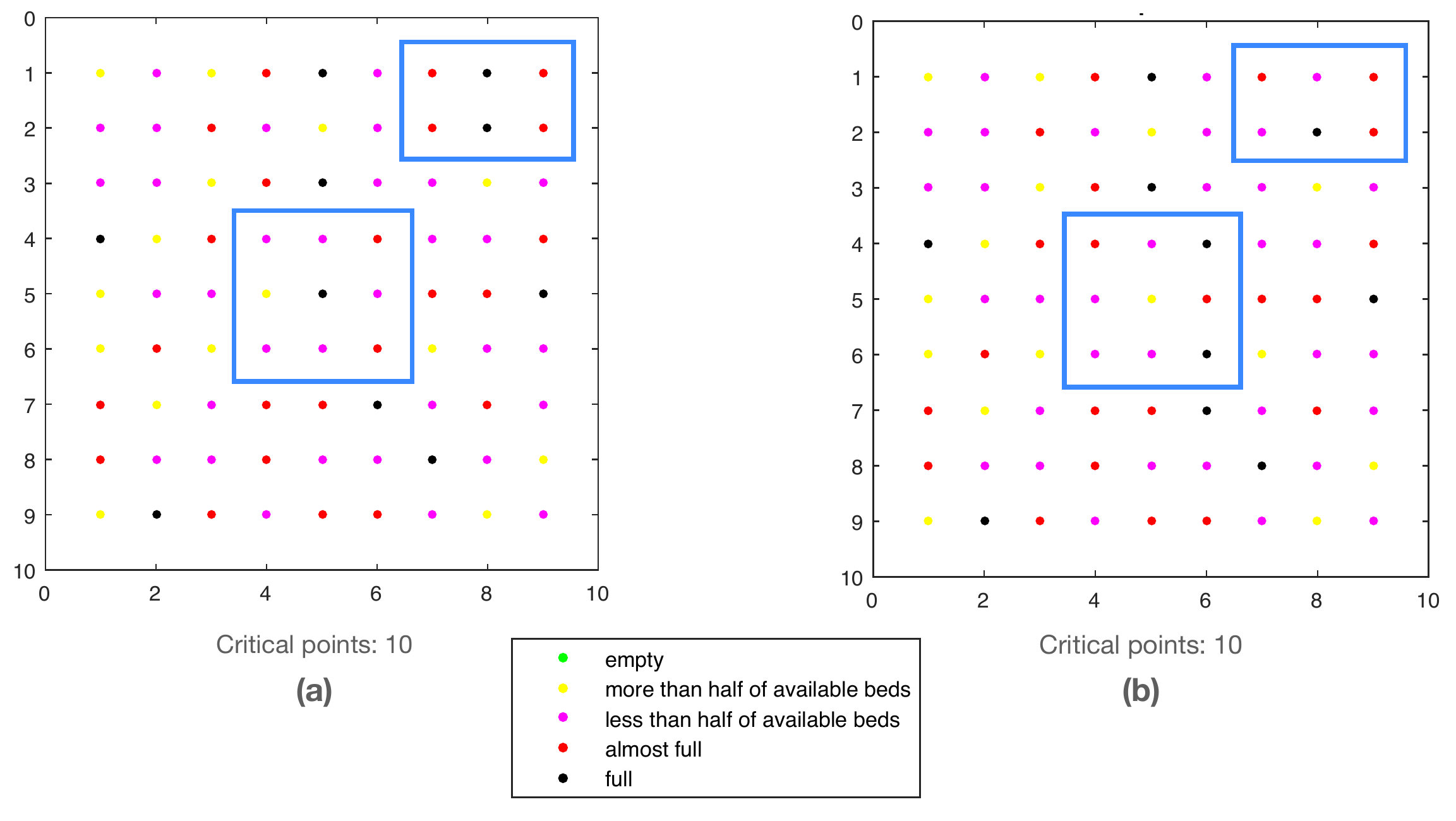}
\caption{\small Model comparison from the case of peripheral outbreaks: final configurations $\Psi(\bar{\mathbf{z}}^0)$ (left) and $\Phi(\bar{\mathbf{z}}^0)$ (right). The colors attached to each cell/node represent its corresponding relative capacity (0=green, 1/2=yellow, 3/4=magenta, 4/5=red, 6/7=black)}
\label{figure:figure8}
\end{figure}

Figure~\ref{figure:figure8} provides a representation of the two configurations $\Psi(\bar{\mathbf{z}}^0)$ and $\Phi(\bar{\mathbf{z}}^0\bigr)$, by indicating the level of saturation reached in each hospital.

Finally, the indicator function $\mathcal{F}$ defined at the beginning of Section~\ref{section:examples} takes the following values for the different configurations,
\begin{equation*}
	\mathcal{F}\bigl(\mathbf{w},\bar{\mathbf{z}}^0\bigr)=7.9,\qquad
	\mathcal{F}\bigl(\mathbf{w},\Psi(\bar{\mathbf{z}}^0)\bigr)=5.9,\qquad
	\mathcal{F}\bigl(\mathbf{w},\Phi(\bar{\mathbf{z}}^0\bigr))=4.4\,,
\end{equation*}
and this fact contributes to the advantages of the SHR approach.


\section{Beyond the Abelian Sandpile paradigm}
\label{section:beyond}
\smallskip

The mathematical model introduced in this article lays the foundation for an optimisation of the healthcare system management. A notable feature of the novel proposal is its \emph{scalability} to various levels of description, and also the possibility of improving the experimental simulation by including more realistic situations and different types of medical facilities (emergency rooms, external care points, ...) into an integrated dynamical system.


\subsection{Sandpiles with internal dissipation}
\label{section:dissipation}
\smallskip

In order to incorporate other relevant elements into the SHR model, such as recovery or (unfortunately) death of patients, we postulate an \emph{elimination mechanism} inherent to the system, which is translated in mathematical formalism by considering the presence of some \emph{dissipation} during the evolution.

For instance, we describe the effect of removing particles from the network by a simple subtraction of a (randomly chosen, but possibly measurable) distribution $\boldsymbol{\zeta}=(\zeta_1,\zeta_2,\dots,\zeta_p)$ with the obvious constraint that
\begin{equation*}
	0\leq \zeta_i\leq  \Phi(\mathbf{z}^0)_{i}\qquad \textrm{for any index} \quad i = 1,2,\dots,p\,.
\end{equation*}
The modified model workflow is divided into the following steps.
\begin{itemize}
\item[1-2.] {\it Initialisation/Inflow}. We proceed as in Section~\ref{section:outgoing}\,;
\item[3-4.] {\it Hub/Additional toppling}. We proceed as in Section~\ref{section:outgoing}\,;
\item[5.] {\it Internal dissipation}.
	We subtract the distribution $\boldsymbol{\zeta}$ to the intermediate\\ configuration $\mathbf{z}^1$,
	leading to the final configuration $\Phi(\bar{\mathbf{z}}^0) = \mathbf{z}^1-\boldsymbol{\zeta}$\,;
\item[6.] {\it Iteration}.
	We reinitialize the ground state with $\mathbf{z}^0$ equal to $\Phi(\bar{\mathbf{z}}^0)$,\\ and we
	restart from 2.
\end{itemize}
We refer to this modified model as \emph{Sandpile with Internal Dissipation} (SID).
\smallskip

In the long run, after a large but finite number of iterations, a balance between the inflow and outflow steps has also to be incorporated in order to guarantee conservation of the total number of patients. In particular, it could be useful to add the hypothesis that
\begin{equation*}
	\sum_{n=1}^{N}\sum_{i=1}^{p} w_i^n = \sum_{n=1}^{N}\sum_{i=1}^{p} \zeta_i^n\,,
\end{equation*}
where $N$ is the total number of iterations, with $\mathbf{w}^n$ and $\boldsymbol{\zeta}^n$ denoting
the inflow and outflow contributions at the time-step $n$, respectively. Indeed, we notice that if the lefthand side is larger than the righthand side, the whole system risks to undergo a finite time collapse, by reaching its \emph{total capability} --sum of the capacities of each single structure-- in finite time.

Other interesting features can be added to provide the model with a higher level of realism: for instance, the presence of some inertia to the transfer process~\cite{HeadRodg97}, giving preference to structures with a certain level of hospitalized patients~\cite{FalkWinkKinz15}, or constrain additional bulk dissipation~\cite{TsucKato00}.


\subsection{Conclusion and perspetives}
\label{section:conclusion}
\smallskip

From the analysis developed in this article, we observe that the standard healthcare system management typically generates highly unstable and unbalanced configurations, where specific geographical areas with semi-empty hospitals alternate with others where all medical facilities are saturated, especially during sudden and unforeseen events like the spread of epidemics. Instead, following the SRH strategy for optimized management of connections between hospitals, seems to produce more sparse allocation of patients, which has to be considered as a preferable configuration in terms of load-balancing.

On the other hand, it is crucial to improve the exchange of information and to provide decision-making tools to the local structures, in order to optimise the healthcare response in normal times and to avoid the collapse of individual hospitals in times of crisis. There are many relevant consequences in the socio-economic field: among others, we stress the riveting possibility of the automation of health protocols, meaning to build an application capable of learning something from the data autonomously, without receiving explicit instructions from the outside.

Such conceptual experimentation could create a learning environment in which policy makers may gain a better understanding of how the system responds to their decisions, providing forecasts of potential
different choices and strategies. We are aware that a paradigm shift is required and we hope to give a contribution to this respect. It is worth stressing that the agreement with realistic experimental data is presently very limited. However, the conceptual framework we have proposed in this article applies in principle to many different contexts, and these research directions are currently being explored.



\small
\section*{Acknowledgements}

The authors thank the Department of Mathematics G. Castelnuovo, Sapienza University of Rome, for hosting the electronic workshop {\it COVID-19 calls for Mathematics} – \url{www1.mat.uniroma1.it/ricerca/ convegni/2020/COVIDchiamaMAT} –\\ which has motivated this article.\\
The authors thank prof. Ferdinando Romano (Department of Public Health and Infectious Diseases, Sapienza University of Rome) for the stimulating discussions about the logistic issues of the healthcare system organisation.\\
The authors also thank prof. Stefano Finzi Vita (Department of Mathematics G. Castelnuovo, Sapienza University of Rome) for useful references on the mathematical theory of the Abelian Sandpile model.





\section*{Authors' contributions}
All authors equally contributed to the conception and design of this article, the acquisition, analysis and interpretation of data, the drafting and revision of the manuscript.


\section*{Consent for publication}
All authors read and approved the submitted manuscript. All authors have agreed to be personally accountable for the contributions, and to ensure that questions related to the accuracy or integrity of any part of this article are appropriately resolved and documented in the literature.









    
\end{document}